\documentclass[11pt]{article}
\usepackage{cite}
\usepackage{mathrsfs} % Ting Kei added
\usepackage{ifthen}
\usepackage{kbordermatrix}
\usepackage{exscale}
\usepackage{tabularx}
\usepackage{latexsym}
\usepackage{amsmath,amssymb,amstext,amsthm} % Lots of math symbols and environments
\usepackage{graphicx,color,epsfig}
\usepackage{graphicx}
\usepackage{fullpage}        % was not available on PC?
\usepackage{float}
\usepackage{verbatim}
\usepackage{seceqn}
\def\R{\mathbb{R}}
\def\H{\mathbb{H}}

\def\eqref#1{{\normalfont(\ref{#1})}}

\def\eqref#1{{\normalfont(\ref{#1})}}
\newcommand{\textdef}[1]{\textit{#1}\index{#1}}
\newcommand{\nul}{\mathrm{null}}
\newcommand{\range}{\mathrm{range}}

\newcommand{\Ss}{{\mathcal S} }

\newcommand{\Mn}{{\mathcal M^n}}

\newcommand{\Mrid}{M_{\real\imag D}}

\newcommand{\A}{{\mathcal A}}

\newcommand{\bbm}{\begin{bmatrix}}
\newcommand{\ebm}{\end{bmatrix}}
\newcommand{\bem}{\begin{pmatrix}}
\newcommand{\eem}{\end{pmatrix}}
\newcommand{\beq}{\begin{equation}}
\newcommand{\beqs}{\begin{equation*}}
\newcommand{\bet}{\begin{table}}
\newcommand{\eeq}{\end{equation}}
\newcommand{\eeqs}{\end{equation*}}
\newcommand{\beqr}{\begin{eqnarray}}

\DeclareMathOperator{\trace}{{trace}}

\DeclareMathOperator{\tr}{{trace}}

\DeclareMathOperator{\Diag}{{Diag}}

\DeclareMathOperator{\sHvec}{{sHvec}}

\DeclareMathOperator{\sHMat}{{sHMat}}

\DeclareMathOperator{\rank}{{rank}}

\DeclareMathOperator{\real}{{\Re}}
\DeclareMathOperator{\imag}{{\Im}}

\newcommand{\nc}{\newcommand}
\nc{\arrow}{{\rm arrow\,}}
\nc{\Arrow}{{\rm Arrow\,}}
\nc{\BoDiag}{{\rm B^0Diag\,}}
\nc{\bodiag}{{\rm b^0diag\,}}

\nc{\Mm}{{\mathcal M}^{m} }
\nc{\Mmn}{{\mathcal M}^{mn} }
\nc{\Mnr}{{\mathcal M}_{nr} }
\nc{\Mnmr}{{\mathcal M}_{(n-1)r} }
\nc{\kwqqp}{Q{$^2$}P\,}
\nc{\kwqqps}{Q{$^2$}Ps}
\def\argmin{\mathop{\rm argmin}}

\nc{\notinaho}{(X,S)\in \overline{AHO}(\A)}
\nc{\inaho}{(X,S)\in AHO(\A)}

\newcommand{\bea}{\begin{eqnarray}}%
\newcommand{\eea}{\end{eqnarray}}%
\newcommand{\beas}{\begin{eqnarray*}}%
\newcommand{\eeas}{\end{eqnarray*}}%
\newcommand{\proj}{{\rm proj}}%
\newcommand{\refl}{{\rm refl}}%
\newcommand{\re}{{\rm Re}\,}%
\newcommand{\im}{{\rm Im}\,}%
{}
\newcommand{\Hnp}[1][]{\,\mathbb{H}_+^{\ifthenelse{\equal{#1}{}}{n}{#1}}}
\newcommand{\Hn}[1][]{\,\mathbb{H}^{\ifthenelse{\equal{#1}{}}{n}{#1}}}
\newcommand{\Dn}[1][]{\,\mathbb{D}^{\ifthenelse{\equal{#1}{}}{n}{#1}}}

\begin{document}
%\begin{multicols}{2}
%\setlength{\multicolsep}{30pt}

%\linenumbers
\bibliographystyle{plain}
\title{
Projection methods in quantum information science
%\footnote{Department
%of Combinatorics and Optimization,
%          Waterloo, Ontario N2L 3G1, Canada,
%         Research Report CORR 2010-01.}
%   \footnote{Research supported by Natural Sciences Engineering Research
%     Council Canada and a grant from AFOSR.}
}
             \author{
%\href{http://orion.math.uwaterloo.ca/~yl2cheun/}
{Vris Yuen-Lam Cheung} \thanks{Department of Combinatorics and Optimization
Faculty of Mathematics, University of Waterloo, Waterloo, Ontario, Canada N2L 3G1; 
%\url{orion.math.uwaterloo.ca/\~yl2cheun}.
}
\and
%\href{http://people.orie.cornell.edu/dd379}
{Dmitriy Drusvyatskiy}\thanks{Department of Mathematics, University of Washington, Seattle, WA 98195-4350; Department of Combinatorics and Optimization, University of Waterloo, Waterloo, Ontario, Canada N2L 3G1; 
%\url{people.orie.cornell.edu/dd379}; 
Research supported by AFOSR.}
\and
%\href{http://people.wm.edu/~cklixx/}
{Chi-Kwong Li}
        \thanks{
Ferguson Professor of Mathematics,
Department of Mathematics, College of William and Mary, Williamsburg, VA 23185;  
%\url{people.wm.edu/\~cklixx}.
}
\and
%\href{https://www.wm.edu/as/appliedscience/graduateprogram/our_students/pelejo_d.php}
{Diane  Christine Pelejo}
        \thanks{
Department of Applied Science, College of William and Mary, Williamsburg, VA 23185; 
%\url{dcpelejo@gmail.com}.
}
\and
%\href{http://orion.math.uwaterloo.ca/~hwolkowi/}
{Henry Wolkowicz}%
        \thanks{Department of Combinatorics and Optimization
        Faculty of Mathematics, University of Waterloo, Waterloo, Ontario, Canada N2L 3G1; Research supported by The Natural Sciences and Engineering Research Council of Canada and by AFOSR; 
%\url{orion.math.uwaterloo.ca/\~hwolkowi}.
}
}
\date{\today}
          \maketitle
%\begin{center}
%          University of Waterloo\\
%          Department of Combinatorics and Optimization\\
%          Waterloo, Ontario N2L 3G1, Canada\\
%          Research Report CORR 2009-04
%\end{center}

%{\bf Key Words:}  

%\noindent {\bf AMS Subject Classification:}

\begin{abstract}
We consider the problem of constructing quantum operations or channels,
if they exist, that
transform a given set of quantum states $\{\rho_1, \dots, \rho_k\}$ to another such set
$\{\hat\rho_1, \dots, \hat\rho_k\}$. In other words,
we must find a {\em completely positive linear map}, if it
exists, that maps a given  set of density matrices to another given set of
density matrices. This problem, in turn, is an instance of a positive semi-definite feasibility problem, but with highly structured constraints. The nature of the constraints makes projection based algorithms very appealing when the number of variables is huge and standard interior point-methods for semi-definite programming are not applicable. We provide emperical evidence to this effect. We moreover present heuristics for finding both high rank and low rank solutions.  
Our experiments are based on the \emph{method of 
alternating projections} and the
\emph{Douglas-Rachford} reflection method.
\end{abstract}

%\linenumbers


\tableofcontents
\listoftables
%\listoffigures

%----------------------------------------------------
\section{Introduction}
%----------------------------------------------------
%In quantum information science, quantum states and quantum operations (a.k.a quantum channels, quantum measurements, etc.) play a fundamental role. 
%Quantum states are mathematically represented as
%\textdef{density matrices} --- positive semidefinite 
%Hermitian matrices with trace one, while
%quantum operations are represented by
%{\em trace preserving completely positive linear maps} on density matrices --- mappings $T$ from the space of $n\times n$ density matrices $\Mn$ to $m\times m$ density matrices $\Mm$ having the form 
%\begin{equation}
%	\label{eq:tpcprepr}
%T(X) = \sum_{j=1}^r F_j X F_j^*,
%\end{equation}
%for some $n\times m$ matrices $F_1, \dots, F_r$ satisfying 
%$\sum_{j=1}^r F^*_j F_j = I_n$. See for example \cite{Choi:75, Kraus:83, NC:00} for more details.

A basic problem in quantum information science is to construct, if it exists,
a \textdef{quantum operation} sending a given set of {\em quantum states} 
$\{\rho_1, \dots, \rho_k\}$ to another set of quantum states
$\{\hat \rho_1, \dots, \hat\rho_k\}$; see e.g.,
\cite{HuangLiPoonSze:12,FungLiSzeChau:13,MR2837768,MR2446039,CJW:04,NC:00} 
and the references therein. Quantum states are mathematically represented as
\textdef{density matrices} --- positive semidefinite 
Hermitian matrices with trace one, while
quantum operations are represented by
{\em trace preserving completely positive linear maps} --- mappings $T$ from the space of $n\times n$ density matrices $\Mn$ to $m\times m$ density matrices $\Mm$ having the form 
\begin{equation}
	\label{eq:tpcprepr}
T(X) = \sum_{j=1}^r F_j X F_j^*,
\end{equation}
for some $n\times m$ matrices $F_1, \dots, F_r$ satisfying 
$\sum_{j=1}^r F^*_j F_j = I_n$. See \cite{Choi:75, Kraus:83, NC:00} for more details.

Thus given some density matrices $A_1, \dots, A_k \in \Mn$ and $B_1, \dots, B_k \in \Mm$, our task is to find a completely positive linear map
$T$ satisfying $T(A_i) = B_i$ for each $i = 1, \dots, k$.
%In mathematical language,
%this problem reduces to finding a \textdef{trace preserving completely positive linear map
%$T$} \footnote{See \cite{Choi:75} for the definition. The exact details will not be important for us, since we will work primarily with the characterization establsihed in \cite{Choi:75}.}, if it
%exists, satisfying $T(A_i) = B_i$ for $i = 1, \dots, k$,
%for given density matrices $A_1, \dots, A_k \in \Mn$
%and $B_1, \dots, B_k \in \Mm$, where
%$\Mn$ denotes the algebra of  $n\times n$ complex matrices.
%By a result of Choi \cite{Choi:75} (see also \cite{Kraus:83}), a 
%\textdef{trace preserving
%completely positive map} $T: \Mn \rightarrow \Mm$
%admits a representation  
%\begin{equation}
%	\label{eq:tpcprepr}
%T(X) = \sum_{j=1}^r F_j X F_j^*,
%\end{equation}
%for some $n\times m$ matrices $F_1, \dots, F_r$ satisfying 
%$\sum_{j=1}^r F_j^* F_j = I_n$.
In turn, if we let $\{E_{11}, E_{12}, \dots, E_{nn}\}$ denote
the \textdef{standard orthonormal basis of $\Mn$}, then
a mapping $T$ is a trace preserving completely positive linear
map if, and only if, the celebrated \textdef{Choi matrix of $T$}, defined in block form by
%the block matrix
\begin{equation}
\label{eq:Choimatrix}
C(T):=\begin{bmatrix}
P_{11}&\ldots &P_{1n}\\
\vdots&P_{ij} & \vdots\\
P_{11}&\ldots &P_{nn}\\
\end{bmatrix}
:= \begin{bmatrix}
T(E_{11})&\ldots &T(E_{1n})\\
\vdots&T(E_{ij}) & \vdots\\
T(E_{11})&\ldots &T(E_{nn})\\
\end{bmatrix} 
\end{equation}
is positive semidefinite and the trace preserving constraints, $\tr(P_{ij}) = \delta_{ij}$, hold, where $\delta_{ij}$ is the 
\textdef{Kronecker delta}. Note that the Choi matrix $C(T)$ is a square $nm\times nm$ matrix, and hence can be very large even for moderate values of $m$ and $n$.
%Here $\Mmn\cong \Mn(\Mm)=\Mm\otimes \Mn$ denotes the
%space of $n\times n$ block-matrices with $m\times m$ matrices as
%blocks. 
%The \textdef{trace preserving conditions} hold if, and
%only if, 
A little thought now shows that our problem is equivalent to the
positive semidefinite feasibility problem for $P=(P_{ij})$:
\begin{equation}
\label{eq:TPCPconstrs}
\begin{array}{cc}
\left\{
\begin{array}{rcl}
\sum_{ij} (A_\ell)_{ij}P_{ij} &=& B_\ell, \quad \ell = 1, \dots, k\\
\tr(P_{ij}) &=& \delta_{ij}, \quad 1\le i \le j \le n  \\
P\in \H^{nm}_+
\end{array}
\right\},%,
\end{array}
\end{equation}
where $\H^{nm}_+$ denotes the space of $nm\times nm$ positive semi-definite Hermitian matrices. 
Moreover, the rank of the Choi matrix $P$ has a natural interpretation: it is equal to the minimal number of summands needed in any representation of the form \eqref{eq:tpcprepr} for the corresponding 
trace preserving completely positive map $T$.

Because of the trace preserving constraints, the solution 
set of \eqref{eq:TPCPconstrs} is bounded. Thus, the problem is never \emph{weakly
infeasible}, i.e.,~infeasible but contains an asymptotically feasible
sequence, e.g.,~\cite{Duff:56}.
In particular, one can use
standard  primal-dual interior point
semidefinite programming packages to  solve the feasibility problem.
However, when the size of the problem ($m,n$) grows, 
the efficiency and especially the accuracy 
of the semidefinite programming approach is limited. To illustrate, even for a reasonable sized problem $m=n=100$, the number of complex variables involved is $10^8/2$. In this paper,
we exploit the special structure of the problem and develop
projection-based methods to solve high dimensional problems
with high accuracy.
We present numerical experiments based on the 
\emph{alternating projection (MAP)} and the
\emph{Douglas-Rachford (DR)} projection/reflection methods.
We see that the DR method significantly outperforms MAP for this problem.
Our numerical results show promise of projection-based approaches for many other types of feasibility
problems arising in quantum information science.

\section{Projection methods for constructing quantum channels}
\subsection{General background on projection methods}
We begin by describing the method of alternating projections (MAP) and the Douglas-Rachford method (DR) in full generality. To this end, consider an Euclidean space ${\bf E}$ with an inner product $\langle\cdot,\cdot\rangle$ and norm $\|\cdot\|$. We are interested in finding a point $x$ lying in the intersection of two closed subsets $A$ and $B$ of ${\bf E}$. For example $A$ may be an affine subspace of Hermitian matrices (over the reals) and $B$ may be the convex cone of positive semi-definite Hermitian matrices (over the reals), as in our basic quantum channel problem \eqref{eq:TPCPconstrs}. Projection based methods then presuppose that given a point $x\in{\bf E}$, finding a point in the nearest-point set
$$\proj_A(x)=\argmin_{a\in A} \{\|x-a\|\}$$ 
is easy, as is finding a point in $\proj_B(x)$. When $A$ and $B$ are convex, the nearest-point sets $\proj_A(x)$ and $\proj_B(x)$ are singletons, of course.
 
Given a current point $a_l\in A$, the method of alternating projections then iterates the following two steps
\begin{align*}
&\textrm{choose}\quad b_{l}	\in \proj_B(a_l)\\
&\textrm{choose}\quad a_{l+1}\in \proj_A(b_{l})
\end{align*}
When $A$ and $B$ are convex and there exists a pair of nearest points of $A$ and $B$, the method always generates iterates converging to such a pair. In particular, when the convex sets $A$ and $B$ intersect, the method converges to some point in the intersection $A\cap B$. Moreover, when the relative interiors of $A$ and $B$ intersect, convergence is R-linear with the rate governed by the cosines of the angles between the vectors $a_{l+1}-b_l$ and $a_l-b_l$. For details, see for example \cite{MR2849884,BB93,siam_rev,conv_alt}. When $A$ and $B$ are not convex, analogous convergence guarantees hold, but only if the method is initialized sufficiently close to the intersection \cite{alt_genproj,alt_man,btheory,altus}. 

The Douglas Rachford algorithm takes a more asymmetric approach. Given a point $x\in{\bf E}$, we define the reflection operator
$$\refl_A(x)=\proj_A(x)+(\proj_A(x)-x).$$
The Douglas Rachford algorithm is then a ``reflect-reflect-average''
method; that is, given a current iterate $x_l\in {\bf E}$, it generates the next iterate by the formula
$$x_{l+1}=\frac{x_l+\refl_A(\refl_B(x_l))}{2}.$$
It is known that for convex instances, the ``projected iterates'' converge  \cite{Lions}. 
The rate of convergence, however, is not well-understood. On the other hand, the method has proven to be extremely effective empirically for many types of problems; see for example \cite{fran,thi,phase}.  

The salient point here is that for MAP and DR to be effective in practice, the nearest point mappings $\proj_A$ and $\proj_B$ must be easy to evaluate. We next observe that for the quantum channel construction problem -- our basic problem -- these mappings are indeed fairly easy to compute (especially the projection onto the affine subspace).

\subsection{Computing projections in the quantum channel construction problem}
In the current work, we always consider the space of Hermitian matrices $\H^{nm}$ as an Euclidean space, that is we regard $\H^{nm}$ as an inner product space over the reals in the obvious way. As usual, we then endow $\H^{nm}$ with the Frobenius norm $\|P\|=\sum_{i,j} (\re P_{i,j})^2+(\im P_{i,j})^2$, where $\re P_{i,j}$ and $\im P_{i,j}$ are the real and the complex parts of $P_{i,j}$, respectively.

Recall that our basic problem is to find a Hermitian matrix $P=(P_{ij})$ satisfying 
\begin{equation}
\label{eq:TPCPconstrs2}
\begin{array}{cc}
\left\{
\begin{array}{rcl}
\sum_{ij} (A_\ell)_{ij}P_{ij} &=& B_\ell, \quad \ell = 1, \dots, k\\
\tr(P_{ij}) &=& \delta_{ij}, \quad 1\le i \le j \le n  \\
P\in\H^{nm}_+
\end{array}
\right\}.%,
\end{array}
\end{equation}
We aim to apply MAP and DR to this formulation. To this end, we first need to introduce some notation to help with the exposition. Define the linear mappings $$\mathcal{L}_A(P):=\Big(\sum_{ij} (A_\ell)_{ij}P_{ij}\Big)_l\quad \textrm{ and }\quad\mathcal{L}_T(P)=\Big(\tr(P_{i,j})\Big)_{i,j},$$
and let $$\mathcal{L}(P)=(\mathcal{L}_A(P),\mathcal{L}_T(P)).$$
Moreover assemble the vectors
$$B=(B_1,\ldots,B_k)\quad \textrm{ and } \quad \Delta=(\delta_{i,j})_{i,j}.$$
Thus we aim to find a matrix $P$ in the intersection of $\H^{nm}_+$ with the affine subspace
$$\mathcal{A}:=\{P:\mathcal{L}(P)=(B ,\Delta)\}.$$
Projecting a Hermitian matrix $P$ onto $\H^{nm}_+$ is standard due to
the Eckart-Young Theorem, \cite{EckartYoung36}. Indeed if $P=U^*\Diag(\lambda_1,\ldots,\lambda_{mn})U$ is an eigenvalue decomposition of $P$, then we have
$$\proj_{\H^{mn}_+}(P)=U^*\Diag(\lambda_1^+,\ldots,\lambda_{mn}^+)U,$$
where for any real number $r$, we set $r^+=\max\{0,r\}$. Thus projecting
a Hermitian matrix onto $\H^{mn}_+$ requires a single eigenvalue
decomposition --- a procedure for which there are many efficient and
well-tested codes (e.g.,~\cite{MR2398876}).

We next describe how to perform the projection onto the affine subspace $\mathcal{A}$, that is how to solve the nearest point problem
$$\min \Big\{\frac{1}{2}\|P-\hat{P}\|^2: \mathcal{L}(\hat{P})=(B ,\Delta)\Big\}.$$
Classically, the solution is 
$$\proj_{\mathcal{A}}(P)=P+\mathcal{L}^{\dag}R,$$
where $\mathcal{L}^{\dag}$ is the Moore-Penrose generalized inverse of $\mathcal{L}$ and $R:=(B,\Delta)-\mathcal{L}(P)$ is the residual.
Finding the Moore-Penrose generalized inverse of a large linear mapping,
like the one we have here, can often be time consuming and error prone.
Luckily, the special structure of the affine constraints in our problem
allow us to find $\mathcal{L}^{\dag}$ both very quickly and very
accurately, so that in all our
experiments the time to compute the projection onto $\mathcal{A}$ is
negligible compared to the computational effort needed to perform the eigenvalue decompositions. We now describe how to compute $\mathcal{L}^{\dag}$ in more detail; full details can be found in the supplementary text\cite{ChDrLiPeWomatrixrepr:14}.

Henceforth, we use the matlab command $\sHvec(A_k)$ to denote a vectorization of the matrix $A_k$.
We now construct the matrix $M\in \R^{k\times m^2}$
by declaring
\begin{equation}
\label{eq:Mmatrix}
M^T= \begin{bmatrix} \sHvec(A_1) & \sHvec(A_2)& \ldots &  \sHvec(A_k)
       \end{bmatrix}.
\end{equation}
We then separate $M$ into three blocks
$$     M=\begin{bmatrix} 
M_{\real} & M_{\imag} & M_D
       \end{bmatrix},$$
where $M_D \in \R^{k\times m}$ has rows formed from 
the diagonals of matrices $A_i$, and
$M_{\real}$ and $M_{\imag}$ have rows formed from 
the real and imaginary parts of $A_i$, respectively, for $i=1,\ldots,k$.
Define now the matrices
\begin{equation}
\label{eq:Mmats}
\begin{array}{rcl}
  \Mrid&:=& 
     \begin{bmatrix}  M_{\real}  & -M_{\imag} &  M_D \end{bmatrix},
\\ N_{\Re\Im D}&:=& \begin{bmatrix} \frac 1{\sqrt 2} \begin{bmatrix} M_{\real}  &
M_{\real} & -M_{\imag}& -M_{\imag} \cr
      -M_{\imag} & M_{\imag}& -M_{\real}  & M_{\real}   \end{bmatrix} 
          \begin{bmatrix} M_D & 0\cr 0 & M_D \end{bmatrix}
        \end{bmatrix}.
\end{array}
\end{equation}
Permuting the rows and columns of $N_{\Re\Im D}$  in a certain way, described in \cite{ChDrLiPeWomatrixrepr:14}, we obtain a matrix denoted by $N_{final}$.
Then $\mathcal{L}$ can be represented in coordinates (i.e. acting on a vectorization of $P$) in a surprisingly simple way, namely as a matrix:

\begin{equation}
L:=\begin{bmatrix}
I_{t(n-1)} \otimes N_{final} & 0 \cr
 0 & \begin{bmatrix} 
  \begin{bmatrix} 
    I_{n-1}\otimes \Mrid
                                    & 0_{k(n-1),n^2} 
                \end{bmatrix}\cr
   \begin{bmatrix}
   e_n \otimes I_{n^2}
   \end{bmatrix}^T
                \end{bmatrix}
\end{bmatrix},
\end{equation}
where $\otimes$ denotes the Kronecker product, and
$t(n-1)$ denotes the triangular number $t(n-1)=\frac{n(n-1)}{2}$.
Let the matrix $(\Mrid)_{null}$ have orthonormal columns that yield a basis for
$\nul(\Mrid)$, i.e.,
\[
	\nul(\Mrid)=\range((\Mrid)_{null}).
\]
%Both $M$ and $\Mrid$ are full row rank. Therefore, 
%is the same as the unique right inverse.
The generalized inverse of the  top-left block is
trivial to find from $N_{final}$. An explicit
expression for the generalized inverse of the bottom right-block can
also be found. Therefore, we get an explicit blocked structure for the
Moore-Penrose generalized inverse of the complete matrix representation.
\begin{equation}
L^\dag 
\vspace{.1in}=
\begin{bmatrix}
I_{t(n-1)} \otimes \mathcal{N}_{final}^\dag & 0 \cr
 0 & \begin{bmatrix} 
    I_{n-1}\otimes  \Mrid^\dag
                                    & e_{n-1}\otimes (\Mrid)_{null} \cr
   e_{n-1}^T \otimes - \Mrid^\dag
     & I_{n^2}-(n-1)(\Mrid)_{null}
                \end{bmatrix}
\end{bmatrix},
\end{equation}
as claimed. Thus $L^\dag $ is easy to construct by simply stacking
various small matrices together in blocks. Moreover, this means that
both expressions $Lp$ and $L^\dagger R$ can be \emph{vectorized} and
evaluated efficiently and accurately.

\section{Numerical experiments}
\label{sect:numexp}
In this section, we numerically illustrate the effectiveness of the projection/reflection methods for solving quantum channel construction problems. %We solve several types
%of feasibility problems: maximum rank; minimum rank; and
%constrained rank.
The large/huge problems were solved on an 
AMD Opteron(tm) Processor 6168, 1900.089 MHz cpu running LINUX.
The smaller problems were solved using an
Optiplex 9020, Intel(R) Core(TM), i7-4770 CPUs, 3.40GHz,3.40 GHz, RAM 16GB
running windows 7.

For simplicity of exposition, in our numerical experiments, we set $n=m$. Moreover, we will impose the unital constraint $T(I_n) = I_n$, a common condition in quantum information science.
%One may impose the assumption that $T(I_n/n) = I_m/m$. In particular,
%if $m = n$, it is known as the unital assumption, i.e.,~$T(I_n) = I_n$.
%Clearly, one only needs to add one matrix to both the 
%$A_1, \dots, A_k$ and to the $B_1, \dots, B_k$.
%In this paper, for simplicity, 
%we assume that both the unital and trace preserving
%constraints are always included. This implies that $m=n$ and 
We note in passing that the unital constraint implies 
that the last constraint in each density matrix block 
of constraints for each $i$ is redundant.
To generate random instances for our tests we proceed as follows.
We start with given integers $m=n,k$ and a value for $r$.
We generate  a Choi matrix $P$ using $r$ random unitary matrices 
$F_i, i=1,\ldots,r$ and a positive probability distribution $d$, i.e.,~we set
\[
P=\sum_{i=1}^rd_iF_iF_i^*.
\]
Note that, given a density matrix $X$, then the 
trace preserving completely positive map can now be
evaluated using the blocked form of $P$ in \eqref{eq:Choimatrix} as
\[
T(X)=\sum_{ij}X_{ij}P_{ij}.
\]
We then generate random density matrices $A_i, i=1,\ldots,k$ and set
$B_i$ as the image of the corresponding 
trace preserving completely positive map $T$
on $A_i$, for all $i$. This
guarantees that we have a feasible instance of rank $r$
and larger/smaller $r$ values result
in larger/smaller rank for the feasible Choi matrix $P$. We set $A_{k+1}$ to
be $I_n$ to enforce the unital constraint.

\subsection{Solving the basic problem with DR}
We first look at our basic feasibility problem \eqref{eq:TPCPconstrs}. We illustrate the numerical results only using the DR algorithm since
we found it to be vastly superior to MAP; see
Section \ref{sect:drmap}, below. We found solutions of huge
problems with surprisingly high accuracy and very few iterations.
The results are presented in Table \ref{table:DRfeasP}.
We give the size of the problem, the number of iterations, the norm of
the residual (accuracy) at the end, the maximum value of the cosine values
indicating the linear rate of convergence, and the total computational time to perform a projection on the PSD cone. 
The projection on the PSD cone dominates the time of the algorithm,
i.e.,~the total time is roughly the number of iterations times the
projection time.
To fathom the size of the problems considered, observe that a problem with $m=n=10^2$ finds a PSD matrix of order $10^4$
which has approximately $10^8/2$ variables. Moreover, 
we reiterate that the solutions are
found with extremely high accuracy in very few iterations.

\begin{table}[h!] 
\centering 
\begin{tabular}{|c|c|c|c|c|} 
\hline 
m=n,k,r & iters & norm-residual & max-cos & PSD-proj-CPUs \\ 
\hline 
90,50,90 & 6   &     5.88e-15     & .7014  &  233.8  \\ 
100,60,90 & 7  &     7.243e-15    & 0.8255     &       821.7 \\
110,65,90 & 7  &     7.983e-15     &  0.8222    &        1484 \\
120,70,90 & 8  &     8.168e-15     &   0.8256 &          2583 \\
130,75,90 & 8  &     7.19e-15    &   0.8288 &          3607 \\
140,80,90 & 9  &     8.606e-15    &   0.8475 &         5832 \\
150,85,90 & 11  &     8.938e-15    &   0.8606 &         6188 \\
160,90,90 & 11  &    9.295e-15   &   0.8718 &         1.079e+04\\
170,95,90 & 12  &    9.412-15   &   0.8918?? &         1.139e+04\\
\hline 
\end{tabular} 
\caption{Using DR algorithm; for solving huge problems} 
\label{table:DRfeasP} 
\end{table} 
Note that the CPU time depends
approximately linearly in the size $m=n$.
%\begin{figure}[htp]
%\centering
%\includegraphics[scale=0.5]{sizeproj.pdf}
%\caption{Size of $m=n$ versus CPU seconds for the PSD projection
%}
%\label{fig:sizevscpu}
%\end{figure}

\subsection{Heuristic for finding max-rank feasible solutions using DR and MAP}
\label{sect:drmap}
We now look at the problem of finding \textdef{high rank feasible solutions}. Recall that this corresponds to finding a trace preserving completely positive map $T$ mapping $A_i$ to $B_i$, so that $T$ necessarily has a long operator sum representation \eqref{eq:tpcprepr}.
We moreover use this section to compare the DR and MAP algorithms.
Our numerical tests fix $m=n,k$ and then change the value of $r$,
i.e.,~the value used to generate the test problems. 

The heuristic for finding a large rank solution starts by finding a (current)
feasible solution $P_c$ using a multiple of
the identity as the starting point $P_0=mnI_{mn}$ and finding a feasible
point $P_c$ using DR. We then set the current point $P_c$ to be the 
barycenter of all the feasible points currently found.
The algorithm then continues by
changing the starting point to the \emph{other side and outside} 
of the PSD cone, i.e.,~the new starting point is found by traveling in 
direction $d=mnI_{mn}-\trace(P_c)P_c$ starting from $P_c$  so that the 
new starting point $P_n:=P_c+\alpha  d$
 is \underline{not} PSD. For instance, we may set
$\alpha=2^i\|d\|^2$ for sufficiently large $i$. 
We then apply the DR algorithm with the new starting point until we find a 
matrix $P\succ 0$ or no increase in the rank occurs.

Again, we see that we find very accurate solutions and solutions of maximum
rank. We find that DR is much more efficient both
in the number of iterations in finding a feasible solution from a given
starting point and in the number of steps in our heuristic needed to
find a large rank solution. In Tables \ref{table:DRmaxitermaxrankP} and \ref{table:MAPmaxitermaxrankP} 
we present the output for
several values of $r$ when using DR and MAP, respectively. 
We use a randomly generated feasibility instance for each value of
$r$ but we start MATLAB with the \emph{rng(default)} settings so the same
random instances are generated.
We note that the DR algorithm is successful for finding a maximum rank
solution and usually after only the first step of the heuristic.
The last three $r=12,10,8$ values required $8,9,12$ steps, respectively.
However, the final $P$ solution was obtained to (a high) $9$ decimal accuracy.
%The cosine of the angle at the end showed that further iterations would
%have resulted in even higher accuracy, i.e. the cosine value was not $1$.

The MAP always requires many more iterations
and at least two steps for the maximum rank solution.
It then fails completely once
$r\leq 12$. In fact, it reaches the maximum number of iterations while
only finding a feasible solution to $3$ decimals accuracy
for $r=12$ and then $2$ decimals accuracy for $r=10,8$.
We see that the cosine value has reached $1$ for $r=12,10,8$ and the
MAP algorithm was making no progress towards convergence.

For each value of $r$ we include:
\begin{enumerate}
	\item
	the number of steps of DR that it took to find the max-rank $P$;
\item
 the minumum/maximum/mean
number of iterations for the steps in finding $P$
\footnote{Note that if the maximum value is the same as $iterlimit$, then the
method failed to attain the desired accuracy $toler$ for this particular
value of $r$.};
\item
	the maximum of the
cosine of the angles between three succesive iterates
\footnote{This  is a good indicator of the expected number of
iterations.};
\item
	the value of the maximum rank found. 
	\footnote{We used the \emph{rank} function in MATLAB with the
default tolerance, i.e.,~$\rank(P)$ is the number of singular values of
$P$ that are larger than  $mn*eps(\|P\|)$, where
$eps(\|P\|)$ is the positive distance from $\|P\|$ to the next larger in
magnitude floating point number of the same precision.
Here we note that we did not 
fail to find a max-rank solution with the DR algorithm.}
\begin{table}[h!] 
\centering 
\begin{tabular}{|c|c|c|c|c|c|c|} 
\hline 
  & rank steps & min-iters & max-iters & mean-iters & max-cos & max rank \\ 
\hline 
r=30 & 1 & 6 & 6 & 6 & 7.008801e-01 & 900 \\ 
\hline 
r=28 & 1 & 7 & 7 & 7 & 7.323953e-01 & 900 \\ 
\hline 
r=26 & 1 & 7 & 7 & 7 & 7.550174e-01 & 900 \\ 
\hline 
r=24 & 1 & 8 & 8 & 8 & 7.911440e-01 & 900 \\ 
\hline 
r=22 & 1 & 9 & 9 & 9 & 8.238539e-01 & 900 \\ 
\hline 
r=20 & 1 & 9 & 9 & 9 & 8.454781e-01 & 900 \\ 
\hline 
r=18 & 1 & 11 & 11 & 11 & 8.730321e-01 & 900 \\ 
\hline 
r=16 & 1 & 15 & 15 & 15 & 8.995266e-01 & 900 \\ 
\hline 
r=14 & 1 & 23 & 23 & 23 & 9.288445e-01 & 900 \\ 
\hline 
r=12 & 8 & 194 & 3500 & 1.916375e+03 & 9.954262e-01 & 900 \\ 
\hline 
r=10 & 9 & 506 & 3500 & 2.605778e+03 & 9.968120e-01 & 900 \\ 
\hline 
r=8 & 12 & 2298 & 3500 & 3.350833e+03 & 9.986002e-01 & 900 \\ 
\hline 
\end{tabular} 
\caption{Using DR algorithm; with $[m~n~k~mn~toler~iterlimit]= [30  ~
30     ~       16    ~       900  ~       1e-14    ~      3500]$;
max/min/mean iter and number rank steps for finding max-rank of $P$. The
$3500$ here means $9$ decimals accuracy attained for last step.} 
\label{table:DRmaxitermaxrankP} 
\end{table} 
\begin{table}[h!] 
\centering 
\begin{tabular}{|c|c|c|c|c|c|c|} 
\hline 
  & rank steps & min-iters & max-iters & mean-iters & max-cos & max rank \\ 
\hline 
r=30 & 2 & 55 & 67 & 61 & 8.233188e-01 & 900 \\ 
\hline 
r=28 & 2 & 65 & 77 & 71 & 8.513481e-01 & 900 \\ 
\hline 
r=26 & 2 & 78 & 89 & 8.350000e+01 & 8.754098e-01 & 900 \\ 
\hline 
r=24 & 2 & 100 & 109 & 1.045000e+02 & 9.040865e-01 & 900 \\ 
\hline 
r=22 & 2 & 124 & 130 & 127 & 9.250665e-01 & 900 \\ 
\hline 
r=20 & 2 & 156 & 158 & 157 & 9.432779e-01 & 900 \\ 
\hline 
r=18 & 2 & 239 & 245 & 242 & 9.689567e-01 & 900 \\ 
\hline 
r=16 & 2 & 388 & 407 & 3.975000e+02 & 9.847052e-01 & 900 \\ 
\hline 
r=14 & 2 & 1294 & 1369 & 1.331500e+03 & 9.980012e-01 & 900 \\ 
\hline 
r=12 & 2 & 3500 & 3500 & \fbox{3500} & 1.000000e+00 & 493 \\ 
\hline 
r=10 & 2 & 3500 & 3500 & \fbox{3500} & 1.000000e+00 & 483 \\ 
\hline 
r=8 & 2 & 3500 & 3500 & \fbox{3500} & 1.000000e+00 & 475 \\ 
\hline 
\end{tabular} 
\caption{Using MAP algorithm; with $[m~ n~ k~ mn~ toler~ iterlimit]= [30
~30   ~          16  ~          900    ~      1e-14  ~         3500]$;
max/min/mean iter and number rank steps for finding max-rank of $P$.
The $3500$ mean-iters means max iterlimit reached; low accuracy
attained.} 
\label{table:MAPmaxitermaxrankP} 
\end{table} 
\end{enumerate}

\subsection{Heuristic for finding low rank and rank constrained solutions}
In quantum information science, one might want to obtain a 
feasible Choi matrix solution 
$P = (P_{ij})$ with low rank, e.g.,~\cite[Section 4.1]{Watrous:0710.0902}.
If we have a bound on the rank, then we
could change the algorithm by adding a rank restriction 
when one projects the current iterate of $P = (P_{ij})$ onto the PSD cone.
That is instead of taking the positive part of $P = (P_{ij})$, we take 
the \emph{nonconvex projection}
\[
P_r:=\sum_{j \le r, \lambda_j > 0} \lambda_j x_jx_j^*,
\]
where $P$ has spectral decomposition $\sum_{j=1}^{mn} \lambda_j x_j x_j^*$ 
with $\lambda_1 \ge \cdots \ge \lambda_{mn}$.

Alternatively, we can do the following.
Suppose a feasible Choi matrix $C(T) = P_c=((P_c)_{ij})$ is found with
$\rank(P_c)=r$.
We can then attempt to find a new Choi matrix of smaller rank restricted
to the face $F$ of the PSD cone where the current $P_c$ is in the
relative interior of $F$, i.e.,~the minimal face of the PSD cone
containing $P_c$.
We do this using facial reduction, e.g.,~\cite{bw1,bw3}.
More specifically, suppose that $P_c=VDV^T$ is a compact spectral
decomposition, where $D\in \Ss^r_{++}$ is diagonal, positive
definite and has rank $r$. Then the minimal face $F$ of the PSD cone containing $P_c$ has the form $F=V\Ss^r_+V^T$.
Recall $Lp=b$ denotes the matrix/vector equation
corresponding to the linear constraints in our basic problem with
$p=\sHvec(P)$. Let $L_{i,:}$ denote the rows of the matrix 
representation $L$. We let $\sHMat=\sHvec^{-1}$. Note that 
$\sHMat=\sHvec^{*}$, the adjoint.
Then each row of the equation $Lp=b$ is equivalent to
\[
\langle L_{i,:}^T,\sHvec(P)\rangle
= \langle \sHMat(L_{i,:}^T),V\bar PV^T\rangle
= \langle V^T\sHMat(L_{i,:}^T)V,\bar P\rangle, \quad \bar P \in \Ss_+^r.
\]
Therefore, we can replace the linear constraints with the smaller system
$\bar L \bar p=b$ with equations $\langle \bar L_{i,:},\bar p\rangle$, where
$\bar L_{i,:}=\sHvec\left(V^T\sHMat(L_{i,:}^T)V\right)$. In addition,
since the current feasible point $P_c$ is in the relative interior of
the face $V\Ss_+^rV^T$, if we start outside the PSD cone $\Ss_+^r$ for
our feasibility search, then we
get a singular feasible $\bar P$ if one exists and so have reduced the
rank of the corresponding initial feasible $P$.
We then repeat this process as long as we get a reduction in the rank.

The MAP approach we are using appears to be especially well suited for
finding low rank solutions. In particular, the facial reduction works
well because we are able to get extremely high accuracy feasible
solutions before applying the compact spectral decomposition.
If the initial $P_0$ that
is projected onto the affine subspace is not positive semidefinite,
then successive iterates on the affine subspace stay outside the
semidefinite cone, i.e.,~we obtain a final feasible solution $\bar P$ that 
is not positive definite if one exists. Therefore, the rank 
of $V\bar V^T$ is reduced from the rank of $P$.
The code for this has been
surprisingly successful in reducing rank. We provide some typical
results for small problems in Table \ref{table:DAMminrank}.
We start with a small rank (denoted by $r$) feasible solution that is used to generate a
feasible problem. Therefore, we know that the minimal rank is $\leq r$.
We then repeatedly solve the problem using facial
reduction until a positive definite solution is found
which means we cannot continue with
the facial reduction. Note that we could restart the algorithm using an
upper bound for the rank obtained from the last rank we obtained.
\begin{table}[h!] 
\centering 
\begin{tabular}{|c|c|c|c|c|} 
\hline 
{m=n,k} & initial rank r & facial red. ranks &  final rank & final norm-residual   \\ 
\hline 
12,10&  11 &     100,50,44,39  &    39        &     1.836e-15   \\
12,10&  10 &     92,61,43,44  &    44        &      1.786e-15 \\
20,14&  20 &     304,105,71  &    71        &      9.648e-15 \\
22,13&  20 &     374,121,75  &    75        &      9.746e-15 \\
\hline
\end{tabular} 
\caption{Using DAM algorithm with facial reduction for decreasing the rank} 
\label{table:DAMminrank} 
\end{table} 

Finally, our tests indicate that the rank constrained problem, which is nonconvex, often can be solved efficiently. Moreover, this problem helps in further
reducing the rank. To see this, suppose that we know a bound, $rbnd$, on the rank of a
feasible $P$. Then, as discussed above,
we change the projection onto the PSD cone by using
only the largest $rbnd$ eigenvalues of $P$. In our tests, if we use $r$,
the value from generating our instances, then we were always successful
in finding a feasible solution of rank $r$. 
Our  final tests appear in Table \ref{table:DAMredrank}.
We generate problems with initial rank $r$. We then start solving a
constrained rank problem with starting constraint rank $r_s$
and decrease this rank by $1$ until we can no longer find a feasible
solution; the final rank with a feasible solution is $r_f$.
\begin{table}[h!] 
\centering 
\begin{tabular}{|c|c|c|c|c|} 
\hline 
${m=n,k}$ & initial rank $r$ & starting constr. rank $r_s$ 
             & final constr. rank $r_f$   \\ 
\hline 
12,9&  15 &     20 &    7           \\
25,16&  35 &     45 &    19           \\
30,21&  38 &     48 &    27           \\
%35,26&  42 &     50 &    ???           \\
\hline
\end{tabular} 
\caption{Using DR algorithm for rank constrained problems with ranks
$r_s$ to $r_f$} 
\label{table:DAMredrank} 
\end{table}

%\cleardoublepage
%\addcontentsline{toc}{section}{Index}
%\label{ind:index}
%\printindex

\addcontentsline{toc}{section}{Bibliography}
%\bibliography{.master,.edm,.psd,.bjorBOOK,.GI,CKbibfile}

\begin{thebibliography}{10}

\bibitem{fran}
A.F.J. Arag\'on, J.M. Borwein, and M.K. Tam.
\newblock Recent results on {D}ouglas-{R}achford methods for combinatorial
  optimization problems.
\newblock {\em Journal of Optimization Theory and Applications}, pages 1--30,
  2013.

\bibitem{BB93}
H.H. Bauschke and J.M. Borwein.
\newblock On the convergence of von {N}eumann's alternating projection
  algorithm for two sets.
\newblock {\em Set-Valued Anal.}, 1(2):185--212, 1993.

\bibitem{siam_rev}
H.H. Bauschke and J.M. Borwein.
\newblock On projection algorithms for solving convex feasibility problems.
\newblock {\em SIAM Rev.}, 38(3):367--426, September 1996.

\bibitem{phase}
H.H. Bauschke, P.L. Combettes, and D.R. Luke.
\newblock Phase retrieval, error reduction algorithm, and {F}ienup variants: a
  view from convex optimization.
\newblock {\em J. Opt. Soc. Amer. A}, 19(7):1334--1345, 2002.

\bibitem{btheory}
H.H. Bauschke, D.R. Luke, H.M. Phan, and X.~Wang.
\newblock Restricted normal cones and the method of alternating projections:
  Theory.
\newblock {\em Set-Valued and Variational Anal.}, pages 1--43, 2013.

\bibitem{bw1}
J.M. Borwein and H.~Wolkowicz.
\newblock Facial reduction for a cone-convex programming problem.
\newblock {\em J. Austral. Math. Soc. Ser. A}, 30(3):369--380, 1980/81.

\bibitem{bw3}
J.M. Borwein and H.~Wolkowicz.
\newblock Regularizing the abstract convex program.
\newblock {\em J. Math. Anal. Appl.}, 83(2):495--530, 1981.

\bibitem{conv_alt}
L.M. Bregman.
\newblock The method of successive projection for finding a common point of
  convex sets.
\newblock {\em Sov. Math. Dokl}, 6:688--692, 1965.

\bibitem{CJW:04}
A.~Chees, R.~Jozsa, and A.~Winter.
\newblock On the existence of physical transformations between sets of quantum
  states.
\newblock {\em Int. J. Quant. Inf.}, 2:11--21, 2004.

\bibitem{ChDrLiPeWomatrixrepr:14}
Y.-L. Cheung, D.~Drusvyatskiy, C.-K. Li, D.C. Pelejo, and H.~Wolkowicz.
\newblock Efficient block matrix representations for the feasible trace
  preserving completely positive problem.
\newblock Technical report, University of Waterloo, Waterloo, Ontario, 2014.
\newblock in progress.

\bibitem{Choi:75}
M.D. Choi.
\newblock Completely positive linear maps on complex matrices.
\newblock {\em Linear Algebra and Appl.}, 10:285--290, 1975.

\bibitem{MR2398876}
J.W. Demmel, O.A. Marques, B.N. Parlett, and C.~V{\"o}mel.
\newblock Performance and accuracy of {LAPACK}'s symmetric tridiagonal
  eigensolvers.
\newblock {\em SIAM J. Sci. Comput.}, 30(3):1508--1526, 2008.

\bibitem{altus}
D.~Drusvyatskiy, A.D. Ioffe, and A.S. Lewis.
\newblock Alternating projections and coupling slope.
\newblock 2014.

\bibitem{Duff:56}
R.J. Duffin.
\newblock Infinite programs.
\newblock In A.W. Tucker, editor, {\em Linear Equalities and Related Systems},
  pages 157--170. Princeton University Press, Princeton, NJ, 1956.

\bibitem{EckartYoung36}
C.~Eckart and G.~Young.
\newblock The approximation of one matrix by another of lower rank.
\newblock {\em Psychometrika}, 1:211--218, 1936.

\bibitem{thi}
V.~Elser, I.~Rankenburg, and P.~Thibault.
\newblock Searching with iterated maps.
\newblock {\em Proceedings of the National Academy of Sciences},
  104(2):418--423, 2007.

\bibitem{MR2849884}
R.~Escalante and M.~Raydan.
\newblock {\em Alternating projection methods}, volume~8 of {\em Fundamentals
  of Algorithms}.
\newblock Society for Industrial and Applied Mathematics (SIAM), Philadelphia,
  PA, 2011.

\bibitem{FungLiSzeChau:13}
C.-H.F. Fung, C.-K. Li, N.-S. Sze, and H.F. Chau.
\newblock Conditions for degradability of tripartite quantum states.
\newblock Technical Report arXiv:1308.6359, University of Hong Kong, 2012.

\bibitem{HuangLiPoonSze:12}
Z.~Huang, C.-K. Li, E.~Poon, and N.-S. Sze.
\newblock Physical transformations between quantum states.
\newblock {\em J. Math. Phys.}, 53(10):102209, 12, 2012.

\bibitem{Kraus:83}
K.~Kraus.
\newblock States, effects, and operations: Fundamental notions of quantum
  theory.
\newblock {\em Lecture Notes in Physics, Springer-Verlag, Berlin}, 190, 1983.

\bibitem{alt_genproj}
A.S. Lewis, D.R. Luke, and J.~Malick.
\newblock Local linear convergence for alternating and averaged nonconvex
  projections.
\newblock {\em Found. Comput. Math.}, 9(4):485--513, 2009.

\bibitem{alt_man}
A.S. Lewis and J.~Malick.
\newblock Alternating projections on manifolds.
\newblock {\em Math. Oper. Res.}, 33(1):216--234, 2008.

\bibitem{MR2837768}
C.-K. Li and Y.-T. Poon.
\newblock Interpolation by completely positive maps.
\newblock {\em Linear Multilinear Algebra}, 59(10):1159--1170, 2011.

\bibitem{MR2446039}
C.-K. Li, Y.-T. Poon, and N.-S. Sze.
\newblock Higher rank numerical ranges and low rank perturbations of quantum
  channels.
\newblock {\em J. Math. Anal. Appl.}, 348(2):843--855, 2008.

\bibitem{Lions}
P.~Lions and B.~Mercier.
\newblock Splitting algorithms for the sum of two nonlinear operators.
\newblock {\em SIAM Journal on Numerical Analysis}, 16(6):964--979, 1979.

\bibitem{NC:00}
M.A. Nielsen and I.L. Chuang, editors.
\newblock {\em Quantum Computation and Quantum Information}.
\newblock Cambridge University Press, 2000.

\bibitem{Watrous:0710.0902}
J.~Watrous.
\newblock Distinguishing quantum operations having few kraus operators.
\newblock {\em Quant. Inf. Comp.}, 8:819--833, 2008.

\end{thebibliography}
\def\cprime{$'$} \def\cprime{$'$} \def\cprime{$'$}
  \def\udot#1{\ifmmode\oalign{$#1$\crcr\hidewidth.\hidewidth
  }\else\oalign{#1\crcr\hidewidth.\hidewidth}\fi} \def\cprime{$'$}
  \def\cprime{$'$} \def\cprime{$'$}
  \def\polhk#1{\setbox0=\hbox{#1}{\ooalign{\hidewidth
  \lower1.5ex\hbox{`}\hidewidth\crcr\unhbox0}}} \def\cprime{$'$}
  \def\polhk#1{\setbox0=\hbox{#1}{\ooalign{\hidewidth
  \lower1.5ex\hbox{`}\hidewidth\crcr\unhbox0}}} \def\cprime{$'$}
  \def\cprime{$'$} \def\cprime{$'$} \def\cprime{$'$} \def\cprime{$'$}
  \def\cprime{$'$} \def\polhk#1{\setbox0=\hbox{#1}{\ooalign{\hidewidth
  \lower1.5ex\hbox{`}\hidewidth\crcr\unhbox0}}} \def\cprime{$'$}
  \def\cprime{$'$} \def\cprime{$'$} \def\cprime{$'$} \def\cprime{$'$}
  \def\cprime{$'$} \def\cprime{$'$} \def\cprime{$'$} \def\cprime{$'$}
  \def\cprime{$'$} \def\cprime{$'$} \def\cprime{$'$} \def\cprime{$'$}
  \def\cprime{$'$} \def\cprime{$'$} \def\cprime{$'$}

\end{document}